\documentclass[a4paper]{article}

\usepackage[utf8]{inputenc}
\usepackage{amscd,amsmath,latexsym,amssymb,amsfonts}
\usepackage{mathrsfs} 
\usepackage{natbib}
\usepackage{verbatim}
\usepackage{graphicx}

\title{Comments on ``Confidence distribution, the frequentist distribution estimator of a parameter --- a
review" by Min-ge Xie and Kesar Singh\footnote{{\em In memoriam:} Between the time I met for the first time
with Prof.~Singh in Rutgers in early April 2012 and the time I wrote this review, Prof.~Singh most sadly passed
away. Although I did not know him well, I think he would have appreciated the intellectual challenge raised in
this intellectual dispute and responded accordingly. I am quite sorry this opportunity will never occur.}}

\author{\sc Christian P. Robert\\
{\em Université Paris-Dauphine, IUF, and CREST, Paris, France}}

\begin{document}

\maketitle

\begin{abstract}
This note is a discussion of the paper ``Confidence distribution" by Min-ge Xie and Kesar Singh, to appear in the
{\em International Statistical Review}.
\end{abstract}

\begin{quote}
{\em ``We have shown how confidence distributions, as a broad concept, can subsume and be associated to
many well-known notions in statistics across different schools of inference."} M. Xie and K. Singh
\end{quote}

I must first acknowledge I am rather baffled about the overall reason of this review and that this bafflement will
necessarily impact the following discussion. Indeed, and this is not truly a coincidence!, I happen (and so
do the authors of the review)  to have
discussed the related paper by \cite{fraser:2011} a few months ago: while I strongly disagreed on the
conclusions of this paper, the central point made by Don Fraser was quite clear, namely to show that Bayesian
posterior statements were ungrounded. The current paper is mostly missing this type of clear message and it does not
convey a true sense of support for using confidence distributions. I find instead that the paper meanders
rather aimlessly around the definition of confidence distributions, which are in short dual representations of
frequentist confidence sets, and that it never reaches any definitive conclusion about the appeal of relying on
those confidence distributions... For instance, I had to wait till Section 4 to be introduced to inference
based on confidence distributions and I find the description anticlimactic: using confidence distributions to 
\begin{itemize}
\item[--] construct confidence intervals is hardly surprising, since this is how those distributions are constructed; 
\item[--] derive point estimators does not show any advance beyond convergence in probability; 
\item[--] conduct testing of hypotheses simply provides a recovery of the usual $p$-value both in the 
one-sided and two-sided cases, and again is hardly surprising given the duality between tests and confidence procedures. 
\end{itemize} 
Similarly, optimality is defined as to mirror UMPU test optimality \citep{lehmann:1986}. The most
fruitful connection witness applications of confidence distributions appear later in Section 7, in particular
with the reinterpretation of bootstrap (Section 7.3), although I am far from convinced that ``the concept of
confidence distribution is much broader" than the one of bootstrap distribution.  Furthermore, the very concept
of confidence distributions is restricted (at least in the paper) to unidimensional entities, and seems to possess as
many avatars as there are ways of constructing confidence intervals. So, by the end of this long review, I do
remain skeptical about the innovation (for frequentist theory) brought by adopting the perspective of
confidence distributions.

\begin{quote}
{\em ``Clearly a confidence distribution does not have to be a Bayes posterior distribution, as there are numerous
ways to derive it."} M. Xie and K. Singh
\end{quote}

The fundamental difficulty I have with confidence distributions is the same I have with fiducial distributions,
namely one of missing a proper target. Some objective Bayes approaches like matching priors (see, e.g.,
\citealp{robert:2001}) are often criticised for having as sole purpose to mimic a frequentist coverage, hence
questioning the relevance of going the Bayesian way. It seems to me that this methodology of confidence
distributions suffers from the same if reverse drawback: as with Fisher's fiducial distributions, one tries to
produce a posterior distribution without following a Bayesian modelling approach, i.e.~without selecting a
reference prior distribution, hence questioning the relevance of {\em not} going the Bayesian way! As a result,
either the constructed confidence distribution corresponds to a valid Bayesian posterior distribution, in which
case it is highly preferable to conduct the choice and assessment of this prior on a preliminary and open basis
(rather than defaulting to an implicit black-box prior). Or the confidence distribution {\em does not}
correspond to a genuine prior distribution, in which case it is then incoherent in terms of mere probability
theory, thus likely to suffer the same woes as most empirical Bayes approaches (like inefficiency and
over-fitting). Things somehow turn for the worse when the authors consider a ``true" prior distribution
$\pi(\theta)$, which may be a confidence distribution resulting from past experiments, and combine it with the
confidence distribution associated with the current data, shying away from the genuine likelihood: if nothing
else, multiplying two densities of the same random variable together is an impossibility from a probabilistic
perspective. 

The object of a confidence distribution does thus remain a full mystery for me, as I do not see how to use it
with any confidence either as {\em a Bayesian procedure} or as a {\em frequentist procedure}. In the former
perspective, it does not necessarily correspond to a prior distribution and to perceive the confidence
distribution as a way ``to assist the development of objective Bayes approaches" is misguided, in that the
corresponding ``priors" (if any) would then be data-dependent, hence loose the basic coherence of the Bayesian
approach. In the latter perspective, having a probability distribution on a fixed parameter $\theta$ does not
make sense. Except when reinterpreting it as a bootstrap distribution, i.e.~with a randomness endowed by the
observation into $\hat\theta$ rather than from the parameter. I must add that the authors of the review do not
indicate that the methodology has met with widespread use, beyond their own circle.

In Section 6.2, the fact that expert opinion is available to build prior distributions would sound to me as the
most natural way to engage into licit Bayesian activities since the construction of this prior is then
validated by the real world. To replace the exact likelihood with a confidence distribution is a way to shoot
oneself in the foot, by throwing away a coherent and valid scheme for another one wasting some of the
information provided by the data (and contradicting Birnbaum's likelihood principle in addition!). 

\begin{quote}
{\em ``The result is a mathematical coincidence, hinged on the normal assumption."} M. Xie and K. Singh
\end{quote}

The authors seem to consider that having genuine posterior distributions turn into exact confidence
distributions cannot have a deeper explanation than being a freak, i.e.~``a mathematical coincidence". While being a
spectator for this kind of exercise, I would think there are deeper reasons for this agreement, first and
foremost in connection with the Bayesian interpretation of the best unbiased estimator of \cite{pitman:1938}.
Furthermore, the work of \cite{welch:peers:1963} shows that prior distributions can be chosen towards an
agreement with the frequentist coverage.

\begin{quote}
{\em ``We can have multiple confidence distributions for the same parameter under any specific setting."}M. Xie and
K. Singh
\end{quote}

As acknowledged by the authors, the notion of confidence distributions suffers from the same taint of {\em
ad-hocquery} as most frequentist (and empirical Bayes, see \citealp{robert:2001}) procedures, namely that the
confidence distribution can be defined in many ways. Section 5 introduces an ordering of those confidence
distributions but it unfortunately is an incomplete ordering, as most frequentist orderings are, and it is thus
unlikely that two arbitrary solutions can be ordered according to this principle. The strong connection with
UMPU tests---whose own optimality proceeds from an unnatural restriction on testing procedures---reflects this
difficulty.

\begin{figure}[ht]
\begin{minipage}[t]{0.5\linewidth}
\centering
\includegraphics[width=\textwidth]{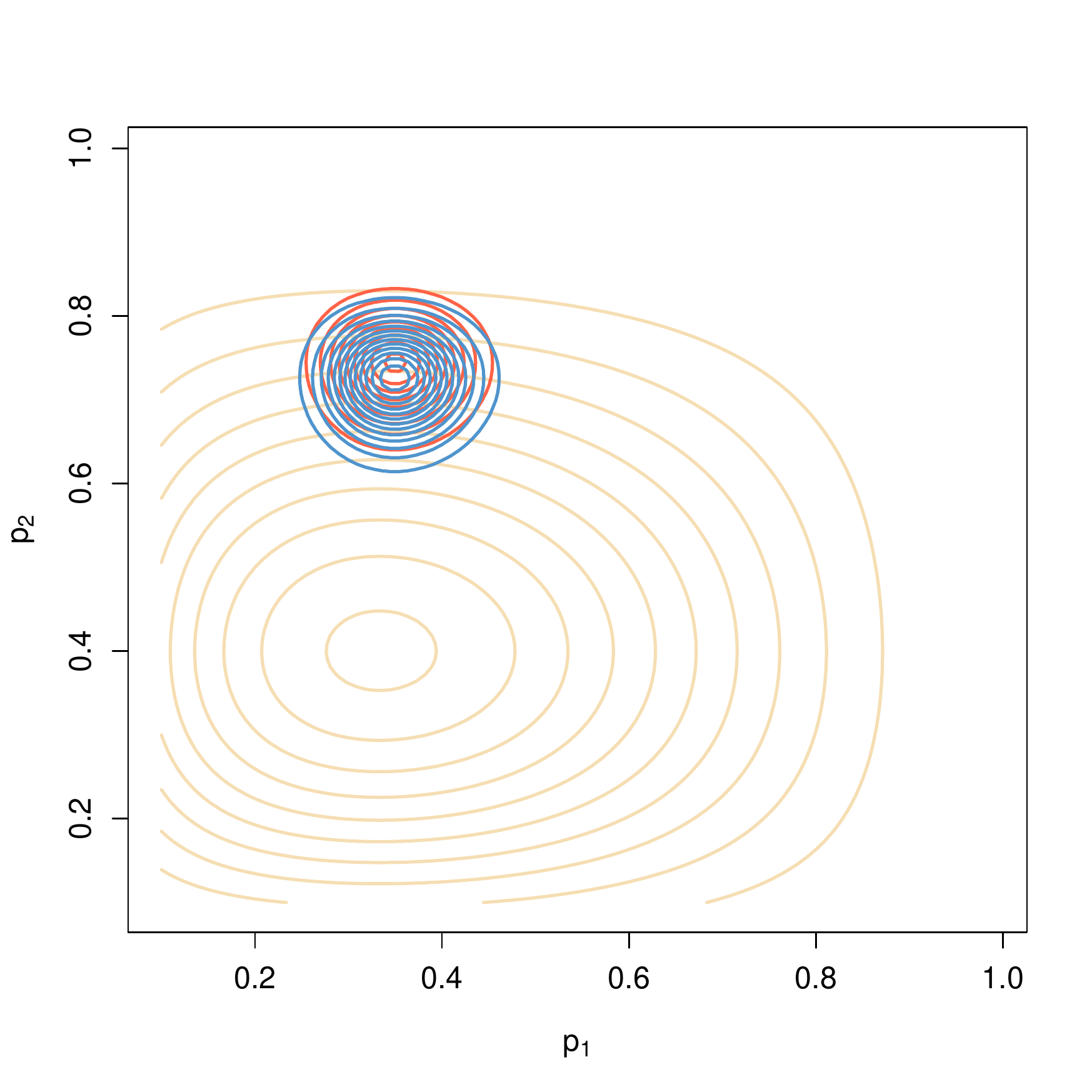}
\caption{Contour plots for the prior (light brown), likelihood (red), and posterior (blue) on $(p_1,p_2)$
based on 100 binomial observations $x_i\sim\mathcal{B}(100,p_i)$ with $x_1=35$ and $x_2=67$ and an
independent prior, $p_1\sim\mathcal{B}e(2,3)$ and $p_2\sim\mathcal{B}e(3,4)$.}
\label{betabin1}
\end{minipage}
\hspace{0.5cm}
\begin{minipage}[t]{0.5\linewidth}
\centering
\includegraphics[width=\textwidth]{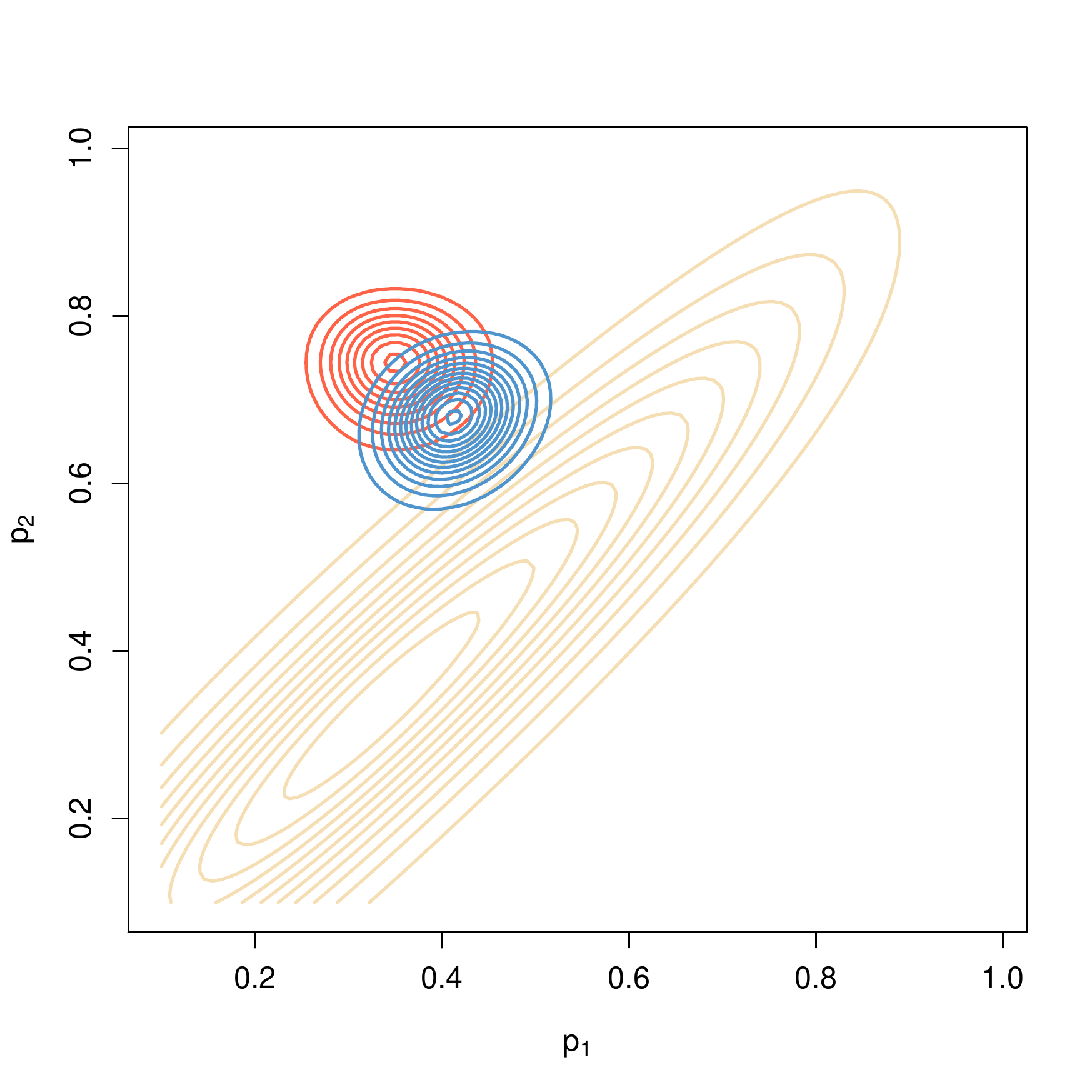}
\caption{Same legend as Figure \ref{betabin1} when using a dependent prior with $p_1\sim\mathcal{B}e(2,3)$ and
$p_2|p_1\sim\mathcal{N}(p_1,0.1)$ restricted to $(0,1)$.}
\label{betabin2}
\end{minipage}
\end{figure}

\begin{quote}
{\em ``A counterintuitive `outlying posterior phenomenon' that is inherent in a Bayesian approach can be avoided in
a confidence distribution-based approach."} M. Xie and K. Singh
\end{quote}

The counter-example discussed in Section 6.2 is only relevant in uncovering the approximation due to the
confidence distribution-based approach, rather than pointing out an inherent flaw in the Bayesian approach.
Indeed, the fact that the posterior distribution is concentrated away from both the prior and the posterior
concentrations seems to be (as far as I can infer given the sparse description contained in the paper) due to
the use of a {\em profile likelihood}, which is an imprecise notion throwing away some of the information
contained in the data. When checking on a regular Bayesian analysis of this beta-binomial model, I could not
spot any discrepancy, using either independent (Fig.~\ref{betabin1}) or dependent (Fig.~\ref{betabin2}) priors.
In any case, the more global issue of having partial prior information like marginal priors on proportions
$p_0$ and $p_1$ does not seem to be such ``a challenging question for Bayesian analysis". Indeed, given those
two marginals, it is always possible to select one parameterised family of copula distributions and to estimate
the parameters of this copula as part of a global Bayesian analysis \citep{silva:lopes:2008}.  

\begin{quote}
{\em ``The review is not intended to re-open the philosophical debate that has lasted more than two hundred
years. On the contrary, it is hoped that the article will help bridge the gap between these different
statistical procedures."} M. Xie and K. Singh
\end{quote}

In conclusion, I fear the authors have not made a proper case in favour of confidence distributions. The notion
carries neither consistency nor optimality features of its own, while it fundamentally relies on the choice of
another frequentist confidence or $p$-value procedure. Worse, the very construction of the confidence
distribution as an inversion of the confidence interval, i.e.~$H_n(\cdot)=\tau_n^{-1}(\cdot)$, reproduces the
common and mislead semantic drift from ``$(-\infty,\tau_n(\alpha)$ contains the true value $\theta_0$ with
probability $\alpha$" to ``$\theta_0$ belongs to the fixed interval $(-\infty,\tau_n(\alpha)$ with probability
$\alpha$".  Further, as reflected by the discussion at the end of Section 6, the review reflects deep
misunderstandings about Bayesian inference. Indeed, speaking of a ``truthful joint prior" or of a ``prior that
is in agreement with the likelihood evidence" shows that the prior is considered as a mythical (if true) unique
entity, rather than as the choice of a reference measure, which is how I do understand priors.


\end{document}